\documentclass{amsart}
\usepackage[all]{xy}

\newtheorem{theorem}{Theorem}[section]
\newtheorem{lemma}{Lemma}[section]
\newtheorem{proposition}{Proposition}[section]

\newtheorem{definition}{Definition}

\theoremstyle{remark}
\newtheorem{remark}{Remark}[section]
\numberwithin{equation}{section}
\numberwithin{lemma}{section}
\numberwithin{proposition}{section}
\numberwithin{definition}{section}

\newcommand{\internalcomment}[1]{}

\begin{document}

\title[On the structure of Calabi-Yau categories]{On the structure of Calabi-Yau categories with a cluster
tilting subcategory}
\author{Gon{\c c}alo Tabuada}
\address{Universit{\'e} Paris 7 - Denis Diderot, UMR 7586
  du CNRS, case 7012, 2 Place Jussieu, 75251 Paris cedex 05, France.}

\thanks{Supported by FCT-Portugal, scholarship {\tt SFRH/BD/14035/2003}.}

\email{
\begin{minipage}[t]{5cm}
tabuada@math.jussieu.fr
\end{minipage}
}

\begin{abstract}
We prove that for $d\geq 2$, an algebraic $d$-Calabi-Yau triangulated
category endowed with a $d$-cluster tilting subcategory is the stable category of
a DG category which is perfectly $(d+1)$-Calabi-Yau and carries a non
degenerate $t$-structure whose heart has enough projectives.
\end{abstract}

\maketitle

\tableofcontents

\section{Introduction}
In this article, we propose a description of a class of
Calabi-Yau categories using the formalism of DG-categories and the
notion of `stabilization', as used for the description of triangulated
orbit categories in section $7$ of \cite{orbit}.
For $d \geq 2$, let $\mathcal{C}$ be an algebraic $d$-Calabi-Yau
triangulated category endowed with a $d$-cluster tilting subcategory
$\mathcal{T}$,
\emph{cf.}~\cite{cluster}~\cite{Iyama32}~\cite{Iyama33}. Such
categories occur for example,
\begin{itemize}
\item[-] in the representation-theoretic approach to
  Fomin-Zelevinsky's cluster algebras \cite{Cluster1},
  \emph{cf.}~\cite{BMRT}~\cite{CK2}~\cite{GLS} and the references
  given there,
\item[-] in the study of Cohen-Macaulay modules over certain isolated
  singularites, \emph{cf.}~\cite{IR}~\cite{cluster}~\cite{IY}, and the
  study of non commutative crepant resolutions~\cite{VdB}, \emph{cf.}~\cite{IR}.
\end{itemize}
From $\mathcal{C}$ and $\mathcal{T}$ we construct an exact dg category
$\mathcal{B}$, which is perfectly $(d+1)$-Calabi-Yau, and a
non-degenerate aisle $\mathcal{U}$,
\emph{cf.}~\cite{t-struture}, in $\mathrm{H}^0(\mathcal{B})$ whose
heart has enough projectives.
We prove, in theorem~\ref{main}, how to recover the category
$\mathcal{C}$ from $\mathcal{B}$ and $\mathcal{U}$ using a
general procedure of stabilization defined in
section~\ref{mainsec}. This extends previous results of
\cite{preprint} to a more general framework. It follows from
\cite{these} that for $d=2$, up to derived equivalence, the category
$\mathcal{B}$ only depends on $\mathcal{C}$ (with its enhancement) and
not on the choice of $\mathcal{T}$.
In the appendix, we show how to naturally extend a $t$-structure, \emph{cf.}~\cite{Ast100}, on the
compact objects of a triangulated category to the whole category.

\section{Acknowledgments}
This article is part of my Ph.~D. thesis under the supervision of
Prof.~B.~Keller. I deeply thank him for countless useful discussions and for his
perfect guidance and generous patience. 

\section{Preliminaries}\label{preli}

Let $k$ be a field. Let $\mathcal{E}$ be a $k$-linear Frobenius
category with split idempotents. Suppose that its stable category
$\mathcal{C}=\underline{\mathcal{E}}$, with suspension functor $S$,
has finite-dimensional $\mathrm{Hom}$-spaces and admits a Serre
functor $\Sigma$, see \cite{B-K}. Let $d\geq 2$ be an
integer. We suppose that $\mathcal{C}$ is Calabi-Yau of CY-dimension
$d$, \emph{i.e.}~\cite{Kontsevich} there is an isomorphism of triangle functors
$$S^d \stackrel{\sim}{\rightarrow} \Sigma\,.$$
We fix such an isomorphism once and for all. See section $4$ of~\cite{cluster} for
several examples of the above situation.

For $X,Y \in \mathcal{C}$
and $n\in \mathbb{Z}$, we put 
$$ \mathrm{Ext}^n(X,Y)=\mathrm{Hom}_{\mathcal{C}}(X,S^nY)\,.$$
We suppose that $\mathcal{C}$ is endowed with a $d$-cluster tilting subcategory $\mathcal{T} \subset \mathcal{C}$, \emph{i.e.}
\begin{itemize}
\item[a)] $\mathcal{T}$ is a $k$-linear subcategory,
\item[b)] $\mathcal{T}$ is functorially finite in $\mathcal{C}$,
  i.e. the functors $\mathrm{Hom}_{\mathcal{C}}(?,X)|\mathcal{T}$ and
  $\mathrm{Hom}_{\mathcal{C}}(X,?)|\mathcal{T}$ are finitely
  generated for all $X \in \mathcal{C}$,
\item[c)] we have $\mathrm{Ext}^i(T,T')=0$ for all $T,T' \in \mathcal{T}$
  and all $0<i<d$ and
\item[d)] if $X \in \mathcal{C}$ satisfies $\mathrm{Ext}^i(T,X)=0$ for all
  $0<i<d$ and all $T \in \mathcal{T}$, then $T$ belongs to $\mathcal{T}$.
 
\end{itemize}
 
Let $\mathcal{M}\subset \mathcal{E}$ be the preimage of $\mathcal{T}$
under the projection functor. In particular, $\mathcal{M}$ contains
the subcategory $\mathcal{P}$ of the projective-injective objects in
$\mathcal{M}$. Note that $\mathcal{T}$ equals the quotient
$\underline{\mathcal{M}}$ of $\mathcal{M}$ by the ideal of morphisms
factoring through a projective-injective.

We dispose of the following commutative square:
$$
\xymatrix{
*+<1pc>{\mathcal{M}} \ar@{^{(}->}[r] \ar@{->>}[d] & \mathcal{E} \ar@{->>}[d] \\
*+<1pc>{\mathcal{T}}  \ar@{^{(}->}[r]  & \underline{\mathcal{E}} =\mathcal{C} \,.
}
$$
We use the notations of \cite{dg-cat-survey}. In particular, for an additive
category $\mathcal{A}$, we denote by $\mathcal{C}(\mathcal{A})$ (resp. $\mathcal{C}^-(\mathcal{A})$, $\mathcal{C}^b(\mathcal{A})$,
$\ldots$) the category of unbounded (resp. right bounded,
resp. bounded, $\ldots$) complexes over $\mathcal{A}$ and by
$\mathcal{H}(\mathcal{A})$ (resp. $\mathcal{H}^-(\mathcal{A})$,
$\mathcal{H}^b(\mathcal{A})$, $\ldots$) its quotient modulo
the ideal of nullhomotopic morphisms. By \cite{KellerVos}, \emph{cf.} also
\cite{Rickard}, the projection functor $\mathcal{E} \rightarrow
\underline{\mathcal{E}}$ extends to a canonical triangle functor
$\mathcal{H}^b(\mathcal{E})/\mathcal{H}^b(\mathcal{P}) \rightarrow
\underline{\mathcal{E}}$. This induces a triangle functor $\mathcal{H}^b(\mathcal{M})/\mathcal{H}^b(\mathcal{P}) \rightarrow
\underline{\mathcal{E}}$. It is shown in \cite{these} that this
functor is a localization functor. Moreover, the projection functor
$\mathcal{H}^b(\mathcal{M}) \rightarrow
\mathcal{H}^b(\mathcal{M})/\mathcal{H}^b(\mathcal{P})$ induces an
equivalence from the subcategory
$\mathcal{H}^b_{\mathcal{E}\mbox{-}ac}(\mathcal{M})$ of bounded
$\mathcal{E}$-acyclic complexes with components in $\mathcal{M}$ onto
its kernel. Thus, we have a short exact sequence of triangulated categories
$$ 0 \longrightarrow \mathcal{H}^b_{\mathcal{E}\mbox{-}ac}(\mathcal{M}) \longrightarrow
\mathcal{H}^b(\mathcal{M})/ \mathcal{H}^b\mathcal(\mathcal{P})
\longrightarrow \mathcal{C} \longrightarrow 0 \,.$$

Let $\mathcal{B}$ be the dg (=differential graded) subcategory of the
category $\mathcal{C}^b(\mathcal{M})_{dg}$ of bounded complexes over
$\mathcal{M}$ whose objets are the $\mathcal{E}$-acyclic complexes. We
denote by $G : \mathcal{H}^-(\mathcal{M}) \rightarrow
\mathcal{D}(\mathcal{B}^{op})^{op}$ the functor which takes a right
bounded complex $X$ over $\mathcal{M}$ to the dg module
$$B \mapsto \mathrm{Hom}^\bullet_{\mathcal{M}}(X,B) \,,$$
where $B$ is in $\mathcal{B}$.
\begin{remark}
By construction, the functor $G$ restricted to
$\mathcal{H}^b_{\mathcal{E}\mbox{-}ac}(\mathcal{M})$ establishes an
equivalence
$$ G : \mathcal{H}^b_{\mathcal{E}\mbox{-}ac}(\mathcal{M})
\stackrel{\sim}{\longrightarrow} \mathrm{per}(\mathcal{B}^{op})^{op} \,.$$

\end{remark}
Recall that if $P$ is a right bounded complex of projectives and $A$
is an acyclic complex, then each morphism from $P$ to $A$ is
nullhomotopic. In particular, the complex
$\mathrm{Hom}^\bullet_{\mathcal{M}}(P,A)$ is nullhomotopic for each $P$
in $\mathcal{H}^- (\mathcal{P})$. Thus $G$ takes $\mathcal{H}^- (P)$ to
zero, and induces a well defined functor (still denoted by $G$)
$$G : \mathcal{H}^b(\mathcal{M})/ \mathcal{H}^b\mathcal(\mathcal{P})
\longrightarrow \mathcal{D}(\mathcal{B}^{op})^{op} \,.$$

\section{Embedding}

\begin{proposition}\label{pleinfidele}
The functor $G$ is fully faithful.
\end{proposition}
For the proof, we need a number of lemmas. \\
It is well-known that the category $\mathcal{H}^-(\mathcal{E})$ admits
a semiorthogonal decomposition, \emph{cf.}~\cite{Bondal-Orlov}, formed by  $\mathcal{H}^-(\mathcal{P})$
and its right orthogonal
$\mathcal{H}^-_{\mathcal{E}\mbox{-}ac}(\mathcal{E})$, the full
subcategory of the right bounded $\mathcal{E}$-acyclic complexes. For $X$ in $\mathcal{H}^-(\mathcal{E})$,
we write 
$$ {\bf{p}}X \rightarrow X \rightarrow {\bf{a}}_pX \rightarrow S{\bf{p}}X$$
for the corresponding triangle, where ${\bf{p}}X$ is in $\mathcal{H}^-(\mathcal{P})$ and
${\bf{a}}_p X$ is in $\mathcal{H}^{-}_{\mathcal{E}\mbox{-}ac}(\mathcal{E})$. 
If $X$ lies in $\mathcal{H}^-(\mathcal{M})$, then clearly
${\bf{a}}_pX$ lies in
$\mathcal{H}^-_{\mathcal{E}\mbox{-}ac}(\mathcal{M})$ so that we have
an induced semiorthogonal decomposition of $\mathcal{H}^-(\mathcal{M})$.

\begin{lemma}\label{fidele}
The functor $\Upsilon:\mathcal{H}^b(\mathcal{M})/ \mathcal{H}^b\mathcal(\mathcal{P})
\longrightarrow \mathcal{H}^-_{\mathcal{E}\mbox{-}ac}(\mathcal{M})$
which takes $X$ to ${\bf{a}}_p X$ is fully faithful.
\end{lemma}

\begin{proof}
By the semiorthogonal decomposition of $\mathcal{H}^-(\mathcal{M})$,
the functor $X \mapsto {\bf{a}}_pX$ induces a right adjoint of the
localization functor
$$\mathcal{H}^-(\mathcal{M}) \longrightarrow \mathcal{H}^-(\mathcal{M})/
\mathcal{H}^-(\mathcal{P})$$
and an equivalence of the quotient category with the right orthogonal $\mathcal{H}^-_{\mathcal{E}\mbox{-}ac}(\mathcal{M})$.
$$
\xymatrix{
 & *+<1pc>{\mathcal{H}^-(\mathcal{P})} \ar@{_{(}->}@<-2ex>[d] & \\
 & *+<1pc>{\mathcal{H}^-(\mathcal{M})} \ar@<-2ex>[d] \ar[u] &
*+<1pc>{\mathcal{H}^{-}_{\mathcal{E}\mbox{-}ac}(\mathcal{M})=\mathcal{H}(\mathcal{P})^{\bot}} \ar@{_{(}->}[l]
\\
*+<1pc>{\mathcal{H}^b(\mathcal{M})/\mathcal{H}^b(\mathcal{P})}  \ar@{^{(}->}[r]
&  *+<1pc>{\mathcal{H}^-(\mathcal{M})/\mathcal{H}^-(\mathcal{P})} \ar@{_{(}->}[u]  \ar@{-->}[ru]^{\sim} & \,.\\
}
$$
Moreover, it is easy
to see that the canonical functor
$$\mathcal{H}^b(\mathcal{M})/\mathcal{H}^b(\mathcal{P}) \longrightarrow
\mathcal{H}^-(\mathcal{M})/\mathcal{H}^-(\mathcal{P})$$
is fully faithful so that we obtain a fully faithful functor
$$ \mathcal{H}^b(\mathcal{M})/\mathcal{H}^b(\mathcal{P}) \longrightarrow
\mathcal{H}^-_{\mathcal{E}\mbox{-}ac}(\mathcal{M})$$
taking $X$ to ${\bf{a}}_p X$.
\end{proof}

\begin{remark}
Since the functor $G$ is triangulated and takes
$\mathcal{H}^-(\mathcal{P})$ to zero, for $X$ in
$\mathcal{H}^b(\mathcal{M})$, the adjunction morphism $X
\rightarrow {\bf{a}}_pX$ yields an isomorphism 
$$ G(X) \stackrel{\sim}{\longrightarrow} G({\bf{a}}_pX)=G(\Upsilon X)\,.$$
\end{remark}

Let $\mathcal{D}^-_{\underline{\mathcal{M}}}(\mathcal{M})$ be the full
subcategory of the derived category $\mathcal{D}(\mathcal{M})$ formed
by the right bounded complexes whose homology modules lie in the
subcategory $\mathrm{Mod}\,\underline{\mathcal{M}}$ of
$\mathrm{Mod}\,\mathcal{M}$. The Yoneda functor $\mathcal{M}
\rightarrow \mathrm{Mod}\,\mathcal{M}$, $M \mapsto M^{\wedge}$, induces a full
imbedding
$$ \Psi :\mathcal{H}^-_{\mathcal{E}\mbox{-}ac}(\mathcal{M}) \hookrightarrow
\mathcal{D}^-_{\underline{\mathcal{M}}}(\mathcal{M}) \,.$$ 
We write $\mathcal{V}$ for its essential image.
Under $\Psi$, the category
$\mathcal{H}^b_{\mathcal{E}\mbox{-}ac}(\mathcal{M})$ is identified
with $\mathrm{per}_{\underline{\mathcal{M}}}(\mathcal{M})$.
Let $\Phi : \mathcal{D}^-_{\underline{\mathcal{M}}}(\mathcal{M}) \rightarrow
\mathcal{D}(\mathcal{B}^{op})^{op}$ be the functor which takes $X$ to the dg module
$$ B \mapsto \mathrm{Hom}^\bullet(X_c, \Psi(B))\,,$$
where $B$ is in $\mathcal{H}^b_{\mathcal{E}\mbox{-}ac}(\mathcal{M})$
and $X_c$ is a cofibrant replacement of $X$ for the projective model
structure on $\mathcal{C}(\mathcal{M})$.
Since for each right bounded complex $M$ with components in
$\mathcal{M}$, the complex $M^{\wedge}$ is cofibrant in
$\mathcal{C}(\mathcal{M})$, it is clear that the functor $G :
\mathcal{H}^b(\mathcal{M})/ \mathcal{H}^b(\mathcal{P}) \rightarrow
\mathcal{D}(\mathcal{B}^{op})^{op}$ is isomorphic to the composition $
\Phi \circ \Psi \circ \Upsilon$.
We dispose of the following commutative diagram
$$
\xymatrix{
\mathcal{H}^b(\mathcal{M})/ \mathcal{H}^b(\mathcal{P})
\ar@{^{(}->}[r]^-{\Upsilon} & \mathcal{H}^-_{\mathcal{E}\mbox{-}ac}(\mathcal{M}) \ar[rr]^{\Psi}
\ar[dr]^{\sim} &  &  \mathcal{D}^-_{\underline{\mathcal{M}}}(\mathcal{M})
\ar[rr]^{\Phi} &  & \mathcal{D}(\mathcal{B}^{op})^{op} \\
 &  & *+<1.5pc>{\mathcal{V}}
  \ar@{^{(}->}[ur] \ar@{_{(}-->}[rrru] & & & \\
 \mathcal{H}^b_{\mathcal{E}\mbox{-}ac}(\mathcal{M}) \ar@{=}[r]
 \ar@{^{(}->}[uu] &  \mathcal{H}^b_{\mathcal{E}\mbox{-}ac}(\mathcal{M}) \ar@{^{(}->}[uu]
\ar[rr]^{\sim} & & *+<1pc>{\mathrm{per}_{\underline{\mathcal{M}}}(\mathcal{M})}
\ar@{^{(}->}[lu] \ar@{^{(}->}[uu] \ar[rr]^{\sim} & &
*+<1pc>{\mathrm{per}(\mathcal{B}^{op})^{op}} \ar@{^{(}->}[uu] \,\,\,.
}
$$
\begin{lemma}\label{caracterisation}
Let $Y$ be an object of $\mathcal{D}^-_{\underline{\mathcal{M}}}(\mathcal{M})$.
\begin{itemize}
\item[a)] $Y$ lies in $\mathrm{per}_{\underline{\mathcal{M}}}(\mathcal{M})$ iff $\mathrm{H}^p(Y)$ is a finitely presented $\underline{\mathcal{M}}$-module
for all $p \in \mathbb{Z}$ and vanishes for all but finitely many $p$.
\item[b)] $Y$ lies in $\mathcal{V}$ iff $\mathrm{H}^p(Y)$ is a finitely presented $\underline{\mathcal{M}}$-module
for all $p \in \mathbb{Z}$ and vanishes for all $p \gg 0$.
\end{itemize}
\end{lemma}

\begin{proof}
a) Clearly the condition is necessary. For the converse, suppose first
that $Y$ is a finitely presented
$\underline{\mathcal{M}}$-module. Then, as an $\mathcal{M}$-module,
$Y$ admits a resolution of length $d+1$ by finitely generated
projective modules by theorem $5.4$ b) of \cite{cluster}. It follows
that $Y$ belongs to $\mathrm{per}_{\underline{\mathcal{M}}}(\mathcal{M})$. Since $\mathrm{per}_{\underline{\mathcal{M}}}(\mathcal{M})$ is
triangulated, it also contains all shifts of finitely presented
$\underline{\mathcal{M}}$-modules and all extensions of shifts. This
proves the converse.\\ 
b) Clearly the condition is necessary.
For the converse, we can suppose
without loss of generality that $Y^n=0$, for all $n \geq 1$ and that
$Y^n$ belongs to $\mathrm{proj}\,\mathcal{M}$, for $n\leq 0$. We now construct a sequence
$$ \cdots \rightarrow P_n \rightarrow \cdots \rightarrow P_1
\rightarrow P_0$$
of complexes of finitely generated projective $\mathcal{M}$-modules
such that $P_n$ is quasi-isomorphic to $\tau_{\geq -n}Y$ for each $n$
and that, for each $p \in \mathbb{Z}$, the sequence of
$\mathcal{M}$-modules $P_n^p$ becomes stationary. By our assumptions,
we have $\tau_{\geq 0}Y \stackrel{\sim}{\rightarrow}
\mathrm{H}^0(Y)$. Since $\mathrm{H}^0(Y)$ belongs to
$\mathrm{mod}\,\underline{\mathcal{M}}$, we know by theorem $5.4$ c) of
\cite{cluster} that it belongs to $\mathrm{per}(\mathcal{M})$ as an
$\mathcal{M}$-module. We define $P_0$ to be a finite resolution of
$\mathrm{H}^0(Y)$ by finitely generated $\mathrm{M}$-modules. 
For the induction step, consider the following truncation triangle associated
with $Y$ 
$$ S^{i+1}\mathrm{H}^{-i-1}(Y) \rightarrow \tau_{\geq -i-1}Y \rightarrow
\tau_{\geq -i}Y \rightarrow S^{i+2}\mathrm{H}^{-i-1}(Y)\,,$$
for $i \geq 0$. 
By the induction hypothesis, we have constructed $P_0,\ldots,P_i$ and we dispose of a
quasi-isomorphism $P_i \stackrel{\sim}{\rightarrow} \tau_{\geq
  -i}Y$. Let $Q_{i+1}$ be a finite resolution of $S^{i+2}\mathrm{H}^{-i-1}(Y)$ by
finitely presented projective $\mathcal{M}$-modules. We dispose of a
morphism $f_i:P_i \rightarrow Q_{i+1}$ and we define
$P_{i+1}$ as the cylinder of $f_i$. We define $P$ as the limit of the
$P_i$ in the category of complexes. We remark that $Y$
is quasi-isomorphic to $P$ and that $P$
belongs to $\mathcal{V}$. This proves the converse.
\end{proof}

Let $X$ be in $\mathcal{H}^-_{\mathcal{E}\mbox{-}ac}(\mathcal{M})$.

\begin{remark}\label{t-strures}
Lemma~\ref{caracterisation} shows that the natural $t$-structure of
$\mathcal{D}(\mathcal{M})$ restricts to a $t$-structure on
$\mathcal{V}$. 
This allows us to express $\Psi(X)$ as
$$ \Psi(X) \stackrel{\sim}{\longrightarrow}
\underset{i}{\mbox{holim}}\,\tau_{\geq-i}\Psi(X) \,,$$
where  $\tau_{\geq-i}\Psi(X)$ is in $\mathrm{per}_{\underline{\mathcal{M}}}(\mathcal{M})$.
\end{remark}

\begin{lemma}\label{holim}
We dispose of the following isomorphism
$$
\Phi(\Psi(X))=\Phi(\underset{i}{\mathrm{holim}}\,\tau_{\geq-i}\Psi(X))
\stackrel{\sim}{\longrightarrow}
\underset{i}{\mathrm{holim}}\,\Phi(\tau_{\geq -i}\Psi(X))\,.
$$
\end{lemma}

\begin{proof}
It is enough to show that the canonical morphism induces a
quasi-isomorphism when evaluated at any object $B$ of
$\mathcal{B}$. We have 
$$\Phi(\underset{i}{\mbox{holim}}\,\tau_{\geq-i}\Psi(X))(B)
=
\mathrm{Hom}^\bullet(\underset{i}{\mbox{holim}}\,\tau_{\geq-i}\Psi(X),
B)\,,$$
but since $B$ is a bounded complex, for each $ n \in \mathbf{Z}$, the sequence
$$ i \mapsto \mathrm{Hom}^n(\tau_{\geq -i}\Psi(X),B) $$
stabilizes as $i$ goes to infinity. This implies that
$$
\mathrm{Hom}^\bullet(\underset{i}{\mbox{holim}}\,\tau_{\geq-i}\Psi(X),B)
\stackrel{\sim}{\longleftarrow} \underset{i}{\mbox{holim}}\,\Phi(\tau_{\geq-i}\Psi(X))(B)\,.
$$
\end{proof}

\begin{lemma}\label{commute}
The functor $\Phi$ restricted to the category $\mathcal{V}$ is fully faithful.
\end{lemma}

\begin{proof}
Let $X,Y$ be in $\mathcal{H}^-_{\mathcal{E}\mbox{-}ac}(\mathcal{M})$. 
The following are canonically isomorphic~:

\begin{eqnarray}
&& \mathrm{Hom}_{\mathcal{D}(\mathcal{B}^{op})^{op}}(\Phi\Psi X,
\Phi\Psi Y)  \nonumber \\
&& \mathrm{Hom}_{\mathcal{D}(\mathcal{B}^{op})}(\Phi\Psi Y,
\Phi\Psi X)  \nonumber \\
&& \mathrm{Hom}_{\mathcal{D}(\mathcal{B}^{op})}(\underset{i}{\mbox{hocolim}}\,\Phi\tau_{\geq-i}\Psi
Y,\underset{j}{\mbox{hocolim}}\,\Phi\tau_{\geq-j}\Psi X
)\\
&& \underset{i}{\mbox{holim}}\,
\mathrm{Hom}_{\mathcal{D}(\mathcal{B}^{op})}(\Phi\tau_{\geq-i}\Psi
Y,\underset{j}{\mbox{hocolim}}\,\Phi\tau_{\geq-j}\Psi X
)  \nonumber \\
&& \underset{i}{\mbox{holim}}\,
  \underset{j}{\mbox{hocolim}}\,
  \mathrm{Hom}_{\mathcal{D}(\mathcal{B}^{op})}(\Phi\tau_{\geq-i}\Psi
  Y,\Phi\tau_{\geq-j}\Psi X) \\
&& \underset{i}{\mbox{holim}}\,
  \underset{j}{\mbox{hocolim}}\,
  \mathrm{Hom}_{\mathrm{per}_{\underline{\mathcal{M}}}(\mathcal{M})}(\tau_{\geq-j}\Psi X,\tau_{\geq-i}\Psi Y) \nonumber \\
&& \underset{i}{\mbox{holim}}\,
  \mathrm{Hom}_{\mathcal{V}}(\underset{j}{\mbox{holim}}\,\tau_{\geq-j}\Psi X,\tau_{\geq-i}\Psi Y) \\
&&
\mathrm{Hom}_{\mathcal{V}}(\Psi(X),\Psi(Y))\,. \nonumber
\end{eqnarray}
Here $(4.1)$ is by the lemma~\ref{holim} seen in
$\mathcal{D}(\mathcal{B}^{op})$, $(4.2)$ is by the fact that
$\Phi\tau_{\geq-i}\Psi Y$ is compact and $(4.3)$ is by the fact that
$\tau_{\geq-i}\Psi Y$ is bounded.
\end{proof}

It is clear now that lemmas~\ref{fidele}, \ref{holim} and \ref{commute} imply the  proposition~\ref{pleinfidele}. 

\section{Determination of the image of $G$}
Let $L_{\rho}:\mathcal{D}^-(\underline{\mathcal{M}}) \rightarrow
\mathcal{D}^-_{\underline{\mathcal{M}}}(\mathcal{M})$ be the
  restriction functor induced by the projection functor $\mathcal{M}
  \rightarrow \underline{\mathcal{M}}$. $L_{\rho}$ admits a left
  adjoint $L: \mathcal{D}^-_{\underline{\mathcal{M}}}(\mathcal{M})
\rightarrow \mathcal{D}^-(\underline{\mathcal{M}})$ which takes $Y$ to
$Y\otimes^{\mathbb{L}}_{\mathcal{M}}\underline{\mathcal{M}}$.
Let $\mathcal{B}^-$ be the dg subcategory of
$\mathcal{C}^-(\mathrm{Mod}\,\mathcal{M})_{dg}$ formed by the objects
of $\mathcal{D}^-_{\underline{\mathcal{M}}}(\mathcal{M})$ that are in the
  essential image of the restriction of $\Psi$ to
$\mathcal{H}^b_{\mathcal{E}\mbox{-}ac}(\mathcal{M})$. Let
$\mathcal{B}'$ be the DG quotient, \emph{cf.}~\cite{Drinfeld}, of
$\mathcal{B}^-$ by its quasi-isomorphisms. It is clear that the dg
categories $\mathcal{B}'$ and $\mathcal{B}$ are quasi-equivalent,
\emph{cf.}~\cite{dg-cat}, and that the natural dg functor
$\mathcal{M} \rightarrow
\mathcal{C}^-(\mathrm{Mod}\,\mathcal{M})_{dg}$ factors through
$\mathcal{B}^-$. Let $R':\mathcal{D}(\mathcal{B}^{op})^{op}
\rightarrow \mathcal{D}(\underline{\mathcal{M}}^{op})^{op}$ be the
restriction functor induced by the dg functor $\underline{\mathcal{M}}
\rightarrow \mathcal{B}'$.
Let $\Phi': \mathcal{D}^-_{\underline{\mathcal{M}}}(\mathcal{M})
\rightarrow \mathcal{D}(\mathcal{B}'^{op})^{op}$ be the functor which
takes $X$ to the dg module 
$$ B' \mapsto \mathrm{Hom}^{\bullet}(X_c,B')\,,$$
where $B'$ is in $\mathcal{B}'$ and $X_c$ is a cofibrant replacement
of $X$ for the projective model structure on
$\mathcal{C}(\mathrm{Mod}\,\mathcal{M})$.
Finally let $\Gamma :\mathcal{D}¯(\underline{\mathcal{M}}) \rightarrow
\mathcal{D}(\underline{\mathcal{M}}^{op})^{op}$ be the functor that
sends $Y$ to 
$$ M \mapsto \mathrm{Hom}^{\bullet}(Y_c,
\underline{\mathcal{M}}(?,M))\,,$$
where $Y_c$ is a cofibrant replacement of $Y$ for the projective model
structure on $\mathcal{C}(\mathrm{Mod}\,\underline{\mathcal{M}})$ and $M$
is in $\underline{\mathcal{M}}$.

We dispose of the following diagram~:
$$
\xymatrix{
& & & \mathcal{D}(\mathcal{B}^{op})^{op} & \mathcal{B} \ar[d] \\
\mathcal{H}^b(\mathcal{M})/\mathcal{H}^b(\mathcal{P})
\ar@{^{(}->}[r]^-{\Upsilon} & \mathcal{H}^-_{\mathcal{E}\mbox{-}ac}(\mathcal{M})
\ar[r]^{\Psi} & \mathcal{D}^-_{\underline{\mathcal{M}}}(\mathcal{M})
\ar[d]_L \ar[r]^{{\Phi}'} \ar[ur]^{\Phi}  & \mathcal{D}(\mathcal{B}'^{op})^{op}
\ar[d]^{R'} \ar[u]_{\sim} & \mathcal{B}' \\
 &  & \mathcal{D}^-(\underline{\mathcal{M}}) \ar[r]_{\Gamma} &
 \mathcal{D}(\underline{\mathcal{M}}^{op})^{op} & \underline{\mathcal{M}} \ar[u] \\
}
$$

\begin{lemma}\label{commutative}
The following square
$$
\xymatrix{
 \mathcal{D}^-_{\underline{\mathcal{M}}}(\mathcal{M}) \ar[r]^{\Phi'}
 \ar[d]_{L} & \mathcal{D}(\mathcal{B}'^{op})^{op} \ar[d]^{R'} &
 \mathcal{B}' \\
\mathcal{D}^-(\underline{\mathcal{M}}) \ar[r]_{\Gamma} &
\mathcal{D}(\underline{\mathcal{M}}^{op})^{op} &
\underline{\mathcal{M}} \ar[u]
}
$$
is commutative.
\end{lemma}

\begin{proof}
By definition $(R' \circ \Phi')(X)(M)$ equals
$\mathrm{Hom}^{\bullet}(X_c,
\underline{\mathcal{M}}(?,\mathcal{M}))$. Since
$\underline{\mathcal{M}}(?,M)$ identifies with $L_{\rho}M^{\wedge}$
and by adjunction, we have 
$$
\mathrm{Hom}^{\bullet}(X_c, \underline{\mathcal{M}}(?,M))
\stackrel{\sim}{\longrightarrow} \mathrm{Hom}^{\bullet}(X_c, L_{\rho}
M^{\wedge}) \stackrel{\sim}{\longrightarrow}
\mathrm{Hom}^{\bullet}((LX)_c, \underline{\mathcal{M}}(?,M))\,,$$
where the last member equals $(\Gamma \circ L)(X)(M)$.
\end{proof}

\begin{lemma}\label{reflet}
The functor $L$ reflects isomorphisms\,.
\end{lemma}

\begin{proof}
Since $L$ is a triangulated functor, it is enough to show that if $L(Y)=0$, then
$Y=0$. Let $Y$ be in
$\mathcal{D}^-_{\underline{\mathcal{M}}}(\mathcal{M})$ such that
$L(Y)=0$. We can suppose, without loss of generality, that
$\mathrm{H}^p(Y)=0$ for all $p>0$. Let us show that
$\mathrm{H}^0(Y)=0$. Indeed, since $\mathrm{H}^0(Y)$ is an
$\underline{\mathcal{M}}$-module, we have $\mathrm{H}^0(Y) \cong
L^0\mathrm{H}^0(Y)$, where $L^0:\mathrm{Mod}\,\mathcal{M} \rightarrow
\mathrm{Mod}\, \underline{\mathcal{M}}$ is the left adjoint of the
inclusion $\mathrm{Mod}\,\underline{\mathcal{M}} \rightarrow
  \mathrm{Mod}\, \mathcal{M}$. Since $\mathrm{H}^p(Y)$ vanishes in
  degrees $p>0$, we have 
$$L^0\mathrm{H}^0(Y) = \mathrm{H}^0(LY)\,.$$
By induction, one concludes that $\mathrm{H}^p(Y)=0$ for all $p \leq 0$.
\end{proof}

\begin{proposition}\label{caracterization}
An objet $Y$ of $\mathcal{D}^-_{\underline{\mathcal{M}}}(\mathcal{M})$
lies in the essential image of the functor $\Psi \circ \Upsilon
: \mathcal{H}^b(\mathcal{M})/ \mathcal{H}^b(\mathcal{P}) \rightarrow
\mathcal{D}^-_{\underline{\mathcal{M}}}(\mathcal{M})$ iff $\tau_{\geq -n}Y$ is in
$\mathrm{per}_{\underline{\mathcal{M}}}(\mathcal{M})$, for all $n \in
\mathbb{Z}$ and $L(Y)$ belongs to $\mathrm{per}(\underline{\mathcal{M}})$.
\end{proposition}

\begin{proof}
Let $X$ be in $\mathcal{H}^b(\mathcal{M})/
\mathcal{H}^b(\mathcal{P})$. By lemma~\ref{caracterisation} a), $\tau_{\geq -n} \Psi \Upsilon (X)$
is in $\mathrm{per}_{\underline{\mathcal{M}}}(\mathcal{M})$, for all
$n \in \mathbb{Z}$. Since $X$ is a bounded complex, there exists an $s
\ll 0$ such that for all $m <s$ the $m$-components of $\Upsilon (X)$ are
in $\mathcal{P}$, which implies that $L \Psi \Upsilon (X)$ belongs to
$\mathrm{per}(\underline{\mathcal{M}})$.\\
Conversely, suppose that $Y$ is an object of
$\mathcal{D}^-_{\underline{\mathcal{M}}}(\mathcal{M})$ which
satisfies the conditions. By lemma~\ref{caracterisation}, $Y$ belongs
to $\mathcal{V}$. Thus we have $Y=\Psi(Y')$ for some $Y'$ in
$\mathcal{H}^-_{\mathcal{E}\mbox{-}ac}(\mathcal{M})$. We now consider
$Y'$ as an object of $\mathcal{H}^-(\mathcal{M})$ and also write
$\Psi$ for the functor $\mathcal{H}^-(\mathcal{M}) \rightarrow
\mathcal{D}^-(\mathcal{M})$ induced by the Yoneda functor.
We can express $Y'$ as 
$$ Y' \stackrel{\sim}{\longleftarrow} \underset{i}{\mbox{hocolim}}\,
\sigma_{\geq -i}Y'\,,$$
where the $\sigma_{\geq -i}$ are the naive truncations. By our
assumptions on $Y'$, $\sigma_{\geq -i}Y'$ belongs to
$\mathcal{H}^b(\mathcal{M})/ \mathcal{H}^b(\mathcal{P})$, for all $i
\in \mathbb{Z}$. The functors $\Psi$ and $L$ clearly commute with the naive
truncations $\sigma_{\geq -i}$ and so we have
$$ L(Y)=L(\Psi Y')  \stackrel{\sim}{\longleftarrow}
\underset{i}{\mbox{hocolim}}\,L(\sigma_{\geq -i}\Psi Y')
\stackrel{\sim}{\longrightarrow}  
\underset{i}{\mbox{hocolim}}\,\sigma_{\geq -i}L(\Psi Y')\,.$$
By our hypotheses, $L(Y)$ belongs to
$\mathrm{per}(\underline{\mathcal{M}})$ and so there exists an $m \gg
0$ such that 
$$ L(Y)= L(\Psi Y') \stackrel{\sim}{\longleftarrow} \sigma_{\geq -m}L(\Psi
Y') = L(\sigma_{\geq -m}\Psi Y')\,.$$
By lemma~\ref{reflet}, the inclusion
$$\Psi(\sigma_{\geq -m}Y)'= \sigma_{\geq -m}\Psi Y' \longrightarrow \Psi(Y')=Y$$
is an isomorphism.
But since $\sigma_{\geq -m}Y'$ belongs to $ \mathcal{H}^b(\mathcal{M})/
\mathcal{H}^b(\mathcal{P})$, $Y$ identifies with $\Psi(\sigma_{\geq -m}Y')$.
\end{proof}

\begin{remark}\label{perfect}
It is clear that if $X$ belongs to
$\mathrm{per}(\underline{\mathcal{M}})$, then $\Gamma(X)$ belongs to
$\mathrm{per}(\underline{\mathcal{M}}^{op})^{op}$. We also have the
following partial converse.
\end{remark}

\begin{lemma}
Let $X$ be in
$\mathcal{D}^-_{\mathrm{mod}\,\underline{\mathcal{M}}}(\underline{\mathcal{M}})$
such that $\Gamma(X)$ belongs to
$\mathrm{per}(\underline{\mathcal{M}}^{op})^{op}$. Then X is in $\mathrm{per}(\underline{\mathcal{M}})$.
\end{lemma}

\begin{proof}
By lemma~\ref{caracterisation} b) we can suppose, without loss of
generality, that $X$ is a right bounded complex with finitely
generated projective components. Applying $\Gamma$, we get a perfect
complex $\Gamma(X)$. In particular $\Gamma(X)$ is homotopic to zero in
high degrees. But since $\Gamma$ is an equivalence 
$$ \mathrm{proj}\, \underline{\mathcal{M}} \stackrel{\sim}{\longrightarrow}
(\mathrm{proj}\,\underline{\mathcal{M}}^{op})^{op}\,,$$
it follows that $X$ is already homotopic to zero in high degrees. 
\end{proof}

\begin{lemma}
The natural left aisle on
$\mathrm{per}_{\underline{\mathcal{M}}}(\mathcal{M})^{op}
\stackrel{\sim}{\rightarrow} \mathrm{per}(\mathcal{B}^{op})$ satisfies
the conditions of proposition~\ref{extension} b).
\end{lemma}

\begin{proof}
Clearly the natural left aisle $\mathcal{U}$ in
$\mathrm{per}_{\underline{\mathcal{M}}}(\mathcal{M})^{op}$ is
non-degenerate. We need to show that for each $C \in
\mathrm{per}_{\underline{\mathcal{M}}}(\mathcal{M})^{op}$, there is an
integer $N$ such that $\mathrm{Hom}(C,S^NU)=0$ for each $U \in
\mathcal{U}$. We dispose of the following isomorphism
$$
\mathrm{Hom}_{\mathrm{per}_{\underline{\mathcal{M}}}(\mathcal{M})^{op}}(C,S^N\mathcal{U})
\stackrel{\sim}{\rightarrow}
\mathrm{Hom}_{\mathrm{per}_{\underline{\mathcal{M}}}(\mathcal{M})}(S^{-N}\mathcal{U}^{op},C)\,,
$$
where $\mathcal{U}^{op}$ denotes the natural right aisle on
$\mathrm{per}_{\underline{\mathcal{M}}}(\mathcal{M})$. Since by
theorem $5.4$ c) of \cite{cluster} an $\underline{\mathcal{M}}$-module
admits a projective resolution of length $d+1$ as an
$\mathcal{M}$-module and $C$ is a bounded complex, we conclude that
for $N \gg 0$
$$
\mathrm{Hom}_{\mathrm{per}_{\underline{\mathcal{M}}}(\mathcal{M})}(S^{-N}\mathcal{U}^{op},C)=0\,.$$
This proves the lemma.
\end{proof}
 
We denote by $\tau_{\leq n}$ and $\tau_{\geq n}$, $n \in \mathbb{Z}$, the
associated truncation functors on $\mathcal{D}(\mathcal{B}^{op})^{op}$.

\begin{lemma}
The functor $\Phi: \mathcal{D}^-_{\underline{\mathcal{M}}}(\mathcal{M}) \rightarrow
\mathcal{D}(\mathcal{B}^{op})^{op}$ restricted to the category
$\mathcal{V}$ is exact with respect to the given $t$-structures.
\end{lemma}

\begin{proof}
We first prove that $\Phi(\mathcal{V}_{\leq 0}) \subset
\mathcal{D}(\mathcal{B}^{op})^{op}_{\leq 0}$. Let $X$ be in
$\mathcal{V}_{\leq 0}$. We need to show that $\Phi(X)$ belongs to
$\mathcal{D}(\mathcal{B}^{op})^{op}_{\leq 0}$. The following have the
same classes of objects~:
\begin{eqnarray}
&& \mathcal{D}(\mathcal{B}^{op})^{op}_{\leq 0} \nonumber \\
&& \mathcal{D}(\mathcal{B}^{op})_{> 0} \nonumber \\
&& (\mathrm{per}(\mathcal{B}^{op})_{\leq 0})^{\bot} \\
&& \overset{\bot}{}(\mathrm{per}(\mathcal{B}^{op})^{op})_{> 0} \,,
\end{eqnarray}
where in $(5.1)$ we consider the right orthogonal in
$\mathcal{D}(\mathcal{B}^{op})$ and in $(5.2)$ we consider the left
orthogonal in $\mathcal{D}(\mathcal{B}^{op})^{op}$. These isomorphisms
show us that $\Phi(X)$ belongs to
$\mathcal{D}(\mathcal{B}^{op})^{op}_{\leq 0}$ iff 
$$ \mathrm{Hom}_{\mathcal{D}(\mathcal{B}^{op})^{op}}(\Phi(X),
\Phi(P))=0\,,$$
for all $P \in \mathrm{per}_{\underline{\mathcal{M}}}(\mathcal{M})_{>
  0}$. 
Now, by lemma~\ref{commute} the functor $\Phi$ is fully faithful and
so
$$ \mathrm{Hom}_{\mathcal{D}(\mathcal{B}^{op})^{op}}(\Phi(X), \Phi(P))
\stackrel{\sim}{\longrightarrow}
\mathrm{Hom}_{\mathrm{per}_{\underline{\mathcal{M}}}(\mathcal{M})}(X,P)\,.$$
Since $X$ belongs to $\mathcal{V}_{\leq 0} $ and $P$ belongs to
$\mathrm{per}_{\underline{\mathcal{M}}}(\mathcal{M})_{> 0}$, we
conclude that
$$\mathrm{Hom}_{\mathrm{per}_{\underline{\mathcal{M}}}(\mathcal{M})}(X,P)=0
\,,$$
which implies that $\Phi(X) \in
\mathcal{D}(\mathcal{B}^{op})^{op}_{\leq 0}$.
Let us now consider $X$ in $\mathcal{V}$. We dispose of the truncation
triangle
$$ \tau_{\leq 0}X \rightarrow X \rightarrow \tau_{>0}X \rightarrow
S\tau_{\leq 0}X \,.$$
The functor $\Phi$ is triangulated and so we dispose of the triangle
$$ \Phi\tau_{\leq 0}X \rightarrow X \rightarrow \Phi\tau_{>0}X \rightarrow
S \Phi \tau_{\leq 0}X \,,$$
where $\Phi\tau_{\leq 0}X$ belonges to
$\mathcal{D}(\mathcal{B}^{op})^{op}_{\leq 0}$.
Since $\Phi$ induces an equivalence between
$\mathrm{per}_{\underline{\mathcal{M}}}(\mathcal{M})$ and
$\mathrm{per}(\mathcal{B}^{op})^{op}$ and $\mathrm{Hom}(P,
\tau_{>0}X) =0$, for all $P$ in $\mathcal{V}_{\leq 0}$, we conclude
that  $\Phi\tau_{>0}X$ belongs to
$\mathcal{D}(\mathcal{B}^{op})^{op}_{> 0}$. This implies the lemma.
\end{proof}

\begin{definition}\label{stable}
Let $\mathcal{D}(\mathcal{B}^{op})^{op}_f$ denote the full triangulated subcategory of $\mathcal{D}(\mathcal{B}^{op})^{op}$
formed by the objects $Y$ such that $\tau_{\geq-n}Y$ is in
$\mathrm{per}(\mathcal{B}^{op})^{op}$, for all $n \in \mathbb{Z}$, and
$R(Y)$ belongs to $\mathrm{per}(\underline{\mathcal{M}}^{op})^{op}$.
\end{definition}

\begin{proposition}\label{caracterisation2}
An objet $Y$ of $\mathcal{D}(\mathcal{B}^{op})^{op}$ lies in the
essential image of the functor $G: \mathcal{H}^b(\mathcal{M})/
\mathcal{H}^b(\mathcal{P}) \rightarrow
\mathcal{D}(\mathcal{B}^{op})^{op}$ iff it belongs to $\mathcal{D}(\mathcal{B}^{op})^{op}_f$.
\end{proposition}

\begin{proof}
Let $X$ be in  $\mathcal{H}^b(\mathcal{M})/
\mathcal{H}^b(\mathcal{P})$. It is clear that the $\tau_{\geq -n} G(X)$ are
in $\mathrm{per}(\mathcal{B}^{op})^{op}$ for all $n \in \mathbb{Z}$. By
proposition~\ref{caracterization} we know that $L \Psi \Upsilon(X)$
belongs to $\mathrm{per}(\underline{\mathcal{M}})$. By lemma~\ref{commutative} and
remark~\ref{perfect} we conclude that $RG(X)$ belongs to $\mathrm{per}(\underline{\mathcal{M}}^{op})^{op}$. 
Let now $Y$ be in $\mathcal{D}(\mathcal{B}^{op})^{op}_f$. We can express
it,  by the dual of lemma~\ref{filtration} as the homotopy limit of the following diagram
$$ \cdots \rightarrow \tau_{\geq-n-1}Y \rightarrow \tau_{\geq-n}Y
\rightarrow \tau_{\geq-n+1}Y \rightarrow \cdots\,,$$
where $\tau_{\geq-n}Y$ belongs to
$\mathrm{per}(\mathcal{B}^{op})^{op}$, for all $n \in \mathbb{Z}$.
But since $\Phi$ induces an equivalence between
$\mathrm{per}_{\underline{\mathcal{M}}}(\mathcal{M})$ and
$\mathrm{per}(\mathcal{B}^{op})^{op}$, this last diagram corresponds
to a diagram
$$ \cdots \rightarrow M_{-n-1} \rightarrow M_{-n} \rightarrow M_{-n+1}
\rightarrow \cdots$$
in $\mathrm{per}_{\underline{\mathcal{M}}}(\mathcal{M})$.
Let $p \in \mathbb{Z}$. The relations among the truncation functors
imply that the image of the above diagram under each homology functor
$\mathrm{H}^p$, $p \in \mathbb{Z}$, is stationary as $n$ goes to $+ \infty$. This implies that 
$$ \mathrm{H}^p\, \underset{n}{\mbox{holim}}\,M_{-n}
\stackrel{\sim}{\longrightarrow}
\underset{n}{\mbox{lim}}\, \mathrm{H}^p\,M_{-n} \cong \mathrm{H}^p\,
M_j\,,
$$
for all $j<p$. We dispose of the following commutative diagram
$$
\xymatrix{
\underset{n}{\mbox{holim}}\,M_{-n} \ar[r] \ar[d] &
\underset{n}{\mbox{holim}}\, \tau_{\geq -i}\,M_{-n} \cong M_{-i}\\
\tau_{\geq -i}\, \underset{n}{\mbox{holim}}\, M_{-n} \ar[ru]^{\sim}
}
$$
which implies that $$ \tau_{\geq -i} \,\underset{n}{\mathrm{holim}} M_{-n}
\stackrel{\sim}{\longrightarrow} M_{-i}\,,$$
for all $i \in \mathbb{Z}$.
Since $\underset{n}{\mbox{holim}}\,M_{-n}$ belongs to $\mathcal{V}$,
lemma~\ref{holim} allows us to conclude that
$\Phi(\underset{n}{\mbox{holim}}\,M_{-n})\cong Y$. We now show that
$\underset{n}{\mbox{holim}}\,M_{-n}$ satisfies the conditions of
proposition~\ref{caracterization}. We know that
$\tau_{\geq -i}\,\underset{n}{\mbox{holim}}\,M_{-n}$ belongs to
$\mathrm{per}_{\underline{\mathcal{M}}}(\mathcal{M})$, for all $i \in \mathbb{Z}$. 
By lemma~\ref{commutative} $(\Gamma \circ
L)(\underset{n}{\mbox{holim}}\,M_{-n})$ identifies with $R(Y)$, which is
in $\mathrm{per}(\underline{\mathcal{M}}^{op})^{op}$.
Since $\underset{n}{\mbox{holim}}\,M_{-n}$ belongs to $\mathcal{V}$,
its homologies lie in $\mathrm{mod}\,\underline{\mathcal{M}}$ and so we
are in the conditions of lemma~\ref{caracterization}, which implies
that $L(\underset{n}{\mbox{holim}}\,M_{-n})$ belongs to
$\mathrm{per}_{\underline{\mathcal{M}}}(\mathcal{M})$. This finishes
the proof.
\end{proof}

\section{Alternative description}
In this section, we present another characterization of the image of
$G$, which was identified as $\mathcal{D}(\mathcal{B}^{op})^{op}_f$
in propositon~\ref{caracterisation2}. Let $M$ denote an object of
$\mathcal{M}$ and also the naturally associated complex in
$\mathcal{H}^b(\mathcal{M})$. Since the category
$\mathcal{H}^b(\mathcal{M})/\mathcal{H}^b(\mathcal{P})$ is generated
by the objects $M \in \mathcal{M}$ and the functor $G$ is fully
faithful, we remark that $\mathcal{D}(\mathcal{B}^{op})^{op}_f$ equals
the triangulated subcategory of $\mathcal{D}(\mathcal{B}^{op})^{op}$
generated by the objects $G(M)$, $M \in \mathcal{M}$. 
The rest of this section is concerned with the problem
of caracterizing the objects $G(M)$, $M \in \mathcal{M}$. We denote by
$P_M$ the projective $\underline{\mathcal{M}}$-module
$\underline{\mathcal{M}}(?,M)$ associated with $M \in \mathcal{M}$ and
by $X_M$ the image of $M$ under $\Psi \circ \Upsilon$.

\begin{lemma}\label{repre}
We dispose of the following isomorphism
$$
\mathrm{Hom}_{\mathcal{D}_{\underline{\mathcal{M}}}^-(\mathcal{M})}(X_M,Y)
\stackrel{\sim}{\longleftarrow}
\mathrm{Hom}_{\mathrm{mod}\,\underline{\mathcal{M}}}(P_M,
\mathrm{H}^0(Y))\,,$$
for all $Y \in \mathcal{D}_{\underline{\mathcal{M}}}^-(\mathcal{M})$.
\end{lemma}

\begin{proof}
Clearly $X_M$ belongs to
$\mathcal{D}_{\underline{\mathcal{M}}}(\mathcal{M})_{\leq 0}$ and is
of the form
$$ \cdots \rightarrow P_n^{\wedge} \rightarrow \cdots \rightarrow P_1^{\wedge} \rightarrow P_0^{\wedge} \rightarrow M^{\wedge} \rightarrow 0\,,$$
where $P_n \in \mathcal{P}$, $n \geq 0$. Now Yoneda's lemma and the
fact that $\mathrm{H}^m(Y)(P_n)=0$, for all $m \in \mathbb{Z}$, $n \geq
0$, imply the lemma.
\end{proof}

\begin{remark}\label{restric}
Since the functor $\Phi$ restricted to $\mathcal{V}$ is fully faithful
and exact, we have
$$ \mathrm{Hom}_{\mathcal{D}(\mathcal{B}^{op})^{op}}(G(M), \Phi(Y))
\stackrel{\sim}{\longleftarrow}
\mathrm{Hom}_{\mathrm{per}(\mathcal{B}^{op})^{op}}(\Phi(P_M),
\mathrm{H}^0(\Phi(Y)))\,,$$
for all $Y \in \mathcal{V}$.
\end{remark}

We now characterize the objects $G(M)=\Phi(X_M)$, $M \in \mathcal{M}$,
in the triangulated category $\mathcal{D}(\mathcal{B}^{op})$. More precisely, we give a description of the functor 
$$ R_M:= \mathrm{Hom}_{\mathcal{D}(\mathcal{B}^{op})}(?,
\Phi(X_M)):\mathcal{D}(\mathcal{B}^{op})^{op} \rightarrow
\mathrm{Mod}\,k$$
using an idea of M.~Van den Bergh, \emph{cf.} lemma~$2.13$ of
\cite{FBR}. Consider the following functor
$$ F_M:=
\mathrm{Hom}_{\mathrm{per}(\mathcal{B}^{op})}(\mathrm{H}^0(?),
\Phi(P_M)): \mathrm{per}(\mathcal{B}^{op})^{op} \rightarrow
\mathrm{mod}\,k\,.$$
 
\begin{remark}
Remark~\ref{restric} shows that the functor $R_M$ when restricted to
$\mathrm{per}(\mathcal{B}^{op})$ coincides with $F_M$.
\end{remark}

Let $DF_M$ be the composition of $F_M$ with the duality
functor $D=\mathrm{Hom}(?,k)$. 
Note that $DF_M$ is homological.

\begin{lemma}\label{DF}
We dispose of the following isomorphism of functors on $\mathrm{per}(\mathcal{B}^{op})$
$$
DF_M \stackrel{\sim}{\longrightarrow}
\mathrm{Hom}_{\mathcal{D}(\mathcal{B}^{op})}(\Phi(X_M), ?[d+1])\,.$$
\end{lemma}

\begin{proof}
The following functors are canonically isomorphic to
$DF \Phi~:$
\begin{eqnarray}
&&
D\mathrm{Hom}_{\mathrm{per}(\mathcal{B}^{op})}(\mathrm{H}^0\Phi(?),
\Phi(P_M)) \nonumber\\
&&
D\mathrm{Hom}_{\mathrm{per}(\mathcal{B}^{op})}(\Phi\mathrm{H}^0(?),
\Phi(P_M)) \\
&&
D\mathrm{Hom}_{\mathrm{per}_{\underline{\mathcal{M}}}(\mathcal{M})}(P_M,
\mathrm{H}^0(?)) \\
&&
D\mathrm{Hom}_{\mathcal{D}_{\underline{\mathcal{M}}}^-(\mathcal{M})}(X_M,?) \\
&&
\mathrm{Hom}_{\mathcal{D}_{\underline{\mathcal{M}}}^-(\mathcal{M})}(?[-d-1],X_M) \\
&&
\mathrm{Hom}_{\mathcal{D}(\mathcal{B}^{op})^{op}}(\Phi(?)[-d-1],\Phi(X_M)) \\
&&
\mathrm{Hom}_{\mathcal{D}(\mathcal{B}^{op})^{op}}(\Phi(X_M),\Phi(?)[d+1]) 
\end{eqnarray} 
 Step $(6.1)$ follows from the fact that $\Phi$ is exact. Step $(6.2)$
 follows from the fact that $\Phi$ is fully faithful and we are
 considering the opposite category. Step $(6.3)$ is a consequence of
 lemma~\ref{repre}. Step $(6.4)$ follows from the $(d+1)$-Calabi-Yau  property and remark~\ref{t-strures}. Step $(6.5)$ is a consequence of $\Phi$ being fully
 faithful and step $(6.6)$ is a consequence of working in the opposite
 category. Since the functor $\Phi^{op}$ establish an equivalence
 between $\mathrm{per}_{\underline{\mathcal{M}}}(\mathcal{M})^{op}$ and
 $\mathrm{per}(\mathcal{B}^{op})$ the lemma is proven.
\end{proof}

Now, we consider the left Kan extension $E_M$ of $DF_M$
along the inclusion $\mathrm{per}(\mathcal{B}^{op}) \hookrightarrow
  \mathcal{D}(\mathcal{B}^{op})$. We dispose of the following
  commutative triangle~:
$$
\xymatrix{
*+<1pc>{\mathrm{per}(\mathcal{B}^{op})} \ar[rr]^{DF_M} \ar@{^{(}->}[d]
& & \mathrm{mod}\,k \\
\mathcal{D}(\mathcal{B}^{op}) \ar@{-->}[urr]_{E_M} & & \,.
}
$$
The functor $E_M$ is homological and preserves coproducts and so
$DE_M$ is cohomological and transforms coproducts into
products. Since $\mathcal{D}(\mathcal{B}^{op})$ is a compactly
generated triangulated category, the Brown representability theorem, \emph{cf.}~\cite{triangular}, implies that there is a $Z_M \in
\mathcal{D}(\mathcal{B}^{op})$ such that 
$$ \mathrm{D} E_M \stackrel{\sim}{\longrightarrow}
\mathrm{Hom}_{\mathcal{D}(\mathcal{B}^{op})}(?, Z_M)\,.$$

\begin{remark}
Since the duality functor $D$ establishes and anti-equivalence
in $\mathrm{mod}\,k$, the functor $DE_M$ restricted to
$\mathrm{per}(\mathcal{B}^{op})$ is isomorphic to $F_M$.
\end{remark} 

\begin{theorem}
We dispose of an isomorphism
$$ G(M) \stackrel{\sim}{\longrightarrow} Z_M\,.$$
\end{theorem}

\begin{proof}
We now construct a morphism of functors from $R_M$ to $DE_M$. Since $R_M$ is representable, by Yoneda's lemma it is enough to
construct an element in $DE_M(\Phi(X_M))$. Let $\mathcal{C}$ be the
category $\mathrm{per}(\mathcal{B}^{op})\downarrow \Phi(X_M)$, whose objects
are the morphisms $Y' \rightarrow \Phi(X_M)$ and let $\mathcal{C}'$ be
the category $X_M \downarrow
\mathrm{per}_{\underline{\mathcal{M}}}(\mathcal{M})$, whose objects
are the morphisms $X_M \rightarrow X'$. The following are
canonically isomorphic~:
\begin{eqnarray}
&& DE_M (\Phi(X_M)) \nonumber \\
&& D \, \underset{\mathcal{C}}{\mbox{colim}}
\mathrm{Hom}_{\mathcal{D}(\mathcal{B}^{op})}(\Phi(X_M), Y'[d+1])\\
&& D \, \underset{\mathcal{C}'}{\mbox{colim}}
\mathrm{Hom}_{\mathcal{D}_{\underline{\mathcal{M}}}^-(\mathcal{M})}(X'[-d-1], X_M)\\
&& \mathrm{D} \,\underset{i}{\mbox{colim}}\, 
\mathrm{Hom}_{\mathcal{D}_{\underline{\mathcal{M}}}^-(\mathcal{M})}((\tau_{\geq
  -i}X_M)[-d-1], X_M) \\
&& \underset{i}{\mbox{lim}} \,\mathrm{D}\mathrm{Hom}_{\mathcal{D}_{\underline{\mathcal{M}}}^-(\mathcal{M})}((\tau_{\geq
  -i}X_M)[-d-1], X_M) \nonumber \\
&& \underset{i}{\mbox{lim}}\,
\mathrm{Hom}_{\mathcal{D}_{\underline{\mathcal{M}}}^-(\mathcal{M})}(X_M,
\tau_{\geq -i}X_M)
\end{eqnarray}

Step $(6.7)$ is a consequence of the definiton of the left Kan
extension and lemma~\ref{DF}. Step $(6.8)$ is obtained by considering
the opposite category. Step $(6.9)$ follows from the fact that the
system $(\tau_{\geq -i}X_M)_{i \in \mathbb{Z}}$ forms a cofinal system
for the index system of the colimit. Step $(6.10)$ follows from the
$(d+1)$-Calabi-Yau property. Now, the image of the identity by the
canonical morphism
$$
\mathrm{Hom}_{\mathcal{D}_{\underline{\mathcal{M}}}^-(\mathcal{M})}(X_M,
X_M) \longrightarrow \underset{i}{\mbox{lim}}\,
\mathrm{Hom}_{\mathcal{D}_{\underline{\mathcal{M}}}^-(\mathcal{M})}(X_M,
\tau_{\geq -i} X_M)\,,$$
give us an element of $(DE_M)(\Phi(X_M))$ and so a morphism
of functors from $R_M$ to $DE_M$. We remark that this
morphism is an isomorphism when evaluated at the objects of
$\mathrm{per}(\mathcal{B}^{op})$. Since both functors $R_M$ and
$DE_M$ are cohomological, transform coproducts into products
and $\mathcal{D}(\mathcal{B}^{op})$ is compactly generated, we
conclude that we dispose of an isomorphism
$$ G(M) \stackrel{\sim}{\longrightarrow} Z_M\,.$$
\end{proof}

\section{The main theorem}\label{mainsec}
Consider the following commutative square as in section~\ref{preli}:
$$
\xymatrix{
*+<1pc>{\mathcal{M}} \ar@{^{(}->}[r] \ar@{->>}[d] & \mathcal{E} \ar@{->>}[d]\\
*+<1pc>{\mathcal{T}}  \ar@{^{(}->}[r]  & \underline{\mathcal{E}}=\mathcal{C}\,.
}
$$
In the previous sections we have constructed, from the above data, a
dg category $\mathcal{B}$ and a left aisle $\mathcal{U}\subset
\mathrm{H}^0(\mathcal{B})$, see~\cite{t-struture}, satisfying the
following conditions~:
\begin{itemize}
\item[-] $\mathcal{B}$ is an exact dg category over $k$ such that
  $\mathrm{H}^0(\mathcal{B})$ has finite-dimensional $\mathrm{Hom}$-spaces and is
  Calabi-Yau of CY-dimension $d+1$,
\item[-] $\mathcal{U} \subset \mathrm{H}^0(\mathcal{B})$ is a
  non-degenerate left aisle such that~:
\begin{itemize}
\item[-] for all $B \in \mathcal{B}$, there is an integer $N$ such
  that $\mathrm{Hom}_{\mathrm{H}^0 (\mathcal{B})}(B,
  S^NU)=0$ for each $U \in \mathcal{U}$,
\item[-] the heart $\mathcal{H}$ of the $t$-structure on
  $\mathrm{H}^0(\mathcal{B})$ associated with
  $\mathcal{U}$ has enough projectives.
\end{itemize}
\end{itemize}

Let now $\mathcal{A}$ be a dg category and $\mathcal{W}\subset
\mathrm{H}^0(\mathcal{A})$ a left aisle satisfying the above
conditions. We can consider the following general construction~:
Let $\mathcal{Q}$ denote the category of projectives of
$\mathcal{H}$. We claim that the following inclusion
$$ \mathcal{Q} \hookrightarrow \mathcal{H} \hookrightarrow
\mathrm{H}^0(\mathcal{A})\,,$$
lifts to a morphism $\mathcal{Q} \stackrel{j}{\rightarrow}
\mathcal{A}$ in the homotopy category of small dg categories,
\emph{cf.}~\cite{dg-cat-survey}~\cite{IMRN}~\cite{cras}.
Indeed, recall the following argument from section~$7$ of
\cite{dg-cat}: Let $\tilde{\mathcal{Q}}$ be the full dg subcategory of $\mathcal{A}$
whose objects are the same as those of $\mathcal{Q}$. Let $\tau_{\leq
  0}\tilde{\mathcal{Q}}$ denote the dg category obtained from
$\tilde{\mathcal{Q}}$ by applying the truncation functor $\tau_{\leq
  0}$ of complexes to each $\mathrm{Hom}$-space. We dispose of the
following diagram in the category of small dg categories
$$
\xymatrix{
 & *+<1pc>{\tilde{\mathcal{Q}}} \ar@{^{(}->}[r] & \mathcal{A} \\
 & \tau_{\leq 0}\tilde{\mathcal{Q}} \ar[u]  \ar[d] & \\
\mathcal{Q} \ar@{=}[r] & \mathrm{H}^0(\tilde{\mathcal{Q}}) & \,.  
}
$$
Let $X$, $Y$ be objects of $\mathcal{Q}$. Since $X$ and $Y$ belong to
the heart of a $t$-structure in $\mathrm{H}^0(\mathcal{A})$, we have
$$ \mathrm{Hom}_{\mathrm{H}^0(\mathcal{A})}(X,Y[-n])=0\,, $$
for $n \geq 1$. The dg category $\mathcal{A}$ is exact, which implies
that 
$$ \mathrm{H}^{-n}
\mathrm{Hom}_{\tilde{\mathcal{Q}}}^{\bullet}(X,Y)
\stackrel{\sim}{\longrightarrow} \mathrm{Hom}_{\mathrm{H}^0(\mathcal{A})}(X,Y[-n])
=0\,,$$
for $n \geq 1$. This shows that the dg functor $\tau_{\leq
  0}\tilde{\mathcal{Q}} \rightarrow \mathrm{H}^0(\tilde{\mathcal{Q}})$
is a quasi-equivalence and so we dispose of a morphism $\mathcal{Q}
\stackrel{j}{\rightarrow} \mathcal{A}$ in the homotopy category of
small dg categories. We dispose of a triangle functor $j^*:\mathcal{D}(\mathcal{A})
\rightarrow \mathcal{D}(\mathcal{Q})$ given by restriction. By
proposition~\ref{extension}, the left aisle $\mathcal{W} \subset
\mathrm{H}^0(\mathcal{A})$ admits a smallest extension to a left aisle
$\mathcal{D}(\mathcal{A}^{op})^{op}_{\leq 0}$ on
$\mathcal{D}(\mathcal{A}^{op})^{op}$.
Let $\mathcal{D}(\mathcal{A}^{op})^{op}_f$ denote the full triangulated subcategory of $\mathcal{D}(\mathcal{A}^{op})^{op}$
formed by the objects $Y$ such that $\tau_{\geq-n}Y$ is in
$\mathrm{per}(\mathcal{A}^{op})^{op}$, for all $n \in \mathbb{Z}$, and
$j^*(Y)$ belongs to $\mathrm{per}(\mathcal{Q}^{op})^{op}$.
 
\begin{definition}
The stable category of $\mathcal{A}$ with respect to $\mathcal{W}$ is
the triangle quotient 
$$\mathrm{stab}(\mathcal{A},\mathcal{W}) = \mathcal{D}(\mathcal{A}^{op})^{op}_f/\mathrm{per}(\mathcal{A}^{op})^{op}\,.$$
\end{definition}

We are now able to formulate the main theorem.
Let $\mathcal{B}$ be the dg category and $\mathcal{U}
\subset\mathrm{H}^0(\mathcal{B})$ the left aisle constructed in
sections $1$ to $5$.
\begin{theorem}\label{main}
The functor $G$ induces an equivalence of categories
$$ \tilde{G}: \mathcal{C} \stackrel{\sim}{\longrightarrow}
\mathrm{stab}(\mathcal{B}, \mathcal{U})\,.$$
\end{theorem}

\begin{proof}
We dispose of the following commutative diagram~:
$$
\xymatrix{
\mathcal{C} \ar@{-->}[rr]^-{\tilde{G}}_-{\sim} & &
\mathrm{stab}(\mathcal{B}, \mathcal{U}) \\
\mathcal{H}^b(\mathcal{M})/\mathcal{H}^b(\mathcal{P}) \ar[u]
\ar[rr]^-G_-{\sim} && \mathcal{D}(\mathcal{B}^{op})^{op}_f \ar[u]\\
\mathcal{H}^b_{\mathcal{E}\mbox{-}ac}(\mathcal{M}) \ar[u]
\ar[rr]_-{\sim} && \mathrm{per}(\mathcal{B}^{op})^{op} \ar[u] \,.\\
}
$$
The functor $G$ is an equivalence since it is fully faithful by
proposition~\ref{pleinfidele} and essentially surjective by
proposition~\ref{caracterisation2}. Since we dispose of an equivalence
$\mathcal{H}^b_{\mathcal{E}\mbox{-}ac}(\mathcal{M})
\stackrel{\sim}{\longrightarrow} \mathrm{per}(\mathcal{B}^{op})^{op}$
by construction of $\mathcal{B}$ and the columns of the above diagram
are short exact sequences of triangulated categories, the theorem is proved.
\end{proof}

\appendix

\section{Extension of $t$-structures}
Let $\mathcal{T}$ be a compactly generated triangulated category with
suspension functor $S$. We denote by
$\mathcal{T}_c$ the full triangulated sub-category of $\mathcal{T}$ formed by the
compact objects, see \cite{triangular}. We use the terminology of \cite{t-struture}. Let
$\mathcal{U} \subseteq \mathcal{T}_c$ be a left aisle.

\begin{proposition}\label{extension}
\begin{itemize}
\item[a)] The left aisle $\mathcal{U}$ admits a smallest extension to a left
aisle  $\mathcal{T}_{\leq 0}$ on $\mathcal{T}$. 
\item[b)] If $\mathcal{U} \subseteq \mathcal{T}_c$ is non-degenerate
  (\emph{i.e.}, $f:X \rightarrow Y$ is invertible iff
  $\mathrm{H}^p(f)$ is invertible for all $p \in \mathbb{Z}$) and for
  each $X \in \mathcal{T}_c$, there is an integer $N$ such that
  $\mathrm{Hom}(X,S^{N}U)=0$ for each $U \in \mathcal{U}$, then
  $\mathcal{T}_{\leq 0}$ is also non-degenerate.
\end{itemize} 
\end{proposition}

\begin{proof}
a) Let $\mathcal{T}_{\leq 0}$ be the smallest full subcategory of
$\mathcal{T}$ that contains $\mathcal{U}$ and is stable under infinite
sums and extensions. It is clear that $\mathcal{T}_{\leq 0}$ is stable
by $S$ since $\mathcal{U}$ is.
We need to show that the inclusion functor $\mathcal{T}_{\leq 0}
\hookrightarrow \mathcal{T}$ admits a right adjoint. For completeness,
we include the following proof, which is a variant of the `small
object argument', \emph{cf.} also~\cite{souto}. We dispose of the following recursive
procedure. Let $X=X_0$ be an object in $\mathcal{T}$. For the initial
step consider all morphisms from any object $P$ in $\mathcal{U}$
to $X_0$. This forms a set $I_0$ since $\mathcal{T}$ is compactly generated
and so we dispose of the following triangle
$$ 
\xymatrix{
\underset{f \in I_0}{\coprod} P \ar[r] &  X_0 \ar[r] &  X_1 \ar@{~>}[r] & \underset{f \in
  I_0}{\coprod} P \,.
}
$$

For the induction step consider the above construction with $X_n$, $n \geq 1$, in the
place of $X_{n-1}$ and $I_n$ in the place of $I_{n-1}$. We dispose of the following diagram
$$
\xymatrix{
X=X_0 \ar[r] & X_1 \ar[r] \ar@{~>}[dl]  & X_2 \ar[r] \ar@{~>}[dl] &
X_3 \ar[r] \ar@{~>}[dl]
& \cdots
\ar[r] & X' \\
\underset{f \in I_0}{\coprod} P
\ar[u] &\underset{f \in I_1}{\coprod} P
\ar[u] & \underset{f \in I_2}{\coprod} P
\ar[u]  &  \underset{f \in I_3}{\coprod} P \ar[u] & & \,,
}
$$
where $X'$ denotes the homotopy colimit of the diagram $(X_i)_{i \in \mathbb{Z}}$. 
Consider now the following triangle
$$ S^{-1}X'  \rightarrow X'' \rightarrow X
\rightarrow X' \,,$$
where the morphism $X \rightarrow X'$ is the transfinite
composition in our diagram.
Let $P$ be in $\mathcal{U}$. We remark that since $P$ is compact,
$\mathrm{Hom}_{\mathcal{T}}(P,X')=0$. This also implies, by
construction of $\mathcal{T}_{\leq 0}$, that
  $\mathrm{Hom}_{\mathcal{T}}(R,X')=0$, for all $R$ in
    $\mathcal{T}_{\leq 0}$.
The long exact sequence obtained by applying the functor
$\mathrm{Hom}_{\mathcal{T}}(R,?)$ to the triangle above shows that
$$ \mathrm{Hom}(R,X'') \stackrel{\sim}{\longrightarrow}
\mathrm{Hom}(R,X)\,.$$
Let $X''_{n-1}$, $n\geq 1$, be an object as in the following triangle
$$ X=X_0 \rightarrow X_n \rightarrow X''_{n-1} \rightarrow S(X)\,.$$
A recursive application of the octahedron axiom implies that
$X''_{n-1}$ belongs to $S(\mathcal{T}_{\leq 0})$, for all $n \geq
1$. We dispose of the isomorphism
$$ \underset{n}{\mbox{hocolim}}\,X''_{n-1}
\stackrel{\sim}{\longrightarrow} S(X'')\,.$$
Since $\underset{n}{\mbox{hocolim}}\,X''_{n-1}$ belongs to
$S(\mathcal{T}_{\leq 0})$, we conclude that $X''$ belongs to
$\mathcal{T}_{\leq 0}$. This shows that the functor that sends $X$ to
$X''$ is the right adjoint of the inclusion functor  $\mathcal{T}_{\leq 0}
\hookrightarrow \mathcal{T}$. 
This proves that $\mathcal{T}_{\leq 0}$ is a left aisle on $\mathcal{T}$.
We now show that the $t$-structure associated to $\mathcal{T}_{\leq
  0}$, \emph{cf.}~\cite{t-struture},  extends, from $\mathcal{T}_c$
to $\mathcal{T}$, the one associated with $\mathcal{U}$. Let $X$ be in $\mathcal{T}_c$. We dispose
of the following truncation triangle associated with $\mathcal{U}$
$$ X_{\mathcal{U}} \rightarrow X \rightarrow X^{\mathcal{U}^{\bot}}
  \rightarrow SX_{\mathcal{U}}\,.$$
Clearly $X_{\mathcal{U}}$ belongs to $\mathcal{T}_{\leq 0}$. We remark
  that $\mathcal{U}^{\bot}=\mathcal{T}_{\leq 0}^{\bot}$, and so  $X^{\mathcal{U}^{\bot}}$ belongs to
  $\mathcal{T}_{>0}:=\mathcal{T}_{\leq 0}^{\bot}$. 

We now show that
  $\mathcal{T}_{\leq 0}$ is the smallest extension of the left aisle $\mathcal{U}$. Let $\mathcal{V}$ be an aisle containing
  $\mathcal{U}$. The inclusion functor $\mathcal{V} \hookrightarrow
  \mathcal{T}$ commutes with sums, because it admits a right
  adjoint. Since $\mathcal{V}$ is stable under extensions and
  suspensions, it contains $\mathcal{T}_{\leq 0}$.\\

b) Let $X$ be in $\mathcal{T}$. We need to show that $X=0$ iff
$\mathrm{H}^p(X)=0$ for all $p \in \mathbb{Z}$. Clearly the condition
is necessary. For the converse, suppose that $\mathrm{H}^p(X)=0$ for
all $p \in \mathbb{Z}$. Let $n$ be an integer. Consider the following
truncation triangle
$$ \mathrm{H}^{n+1}(X) \rightarrow \tau_{> n} X \rightarrow \tau_{>
  n+1}X \rightarrow S\mathrm{H}^{n+1}(X)\,.$$
Since $\mathrm{H}^{n+1}(X)=0$ we conclude that 
$$ \tau_{> n} X \in \underset{m \in \mathbb{Z}}{\bigcap}\mathcal{T}_{>m}\,,$$
for all $n \in \mathbb{Z}$. Now, let $C$ be a compact object of
$\mathcal{T}$. We know that there is a $k \in \mathbb{Z}$ such that $C
\in \mathcal{T}_{\leq k}$. This implies that 
$$ \mathrm{Hom}_{\mathcal{T}}(C, \tau_{> n}X)=0$$ 
for all $n \in \mathbb{Z}$, since $\tau_{> n}X$ belongs to $(\mathcal{T}_{\leq
  k})^{\bot}$. The category $\mathcal{T}$ is compactly generated and so we conclude
  that $\tau_{> n}X=0$, for all $n \in \mathbb{Z}$.
The following truncation triangle
$$ \tau_{\leq n}X \rightarrow X \rightarrow \tau_{> n}X \rightarrow
S\tau_{\leq n}X\,,$$
implies that $\tau_{\leq n} X$ is isomorphic to $X$ for all $n \in
\mathbb{Z}$. This can be rephrased as saying that
$$ X \in \underset{n \in \mathbb{N}}{\bigcap} \mathcal{T}_{\leq
  -n}\,.$$
Now by our hypothesis there is an integer $N$ such that 
$$ \mathrm{Hom}_{\mathcal{T}}(C, \mathcal{U}_{\leq - N})=0\,.$$
Since $C$ is compact and by construction of $\mathcal{T}_{\leq -N}$,
we have
$$ \mathrm{Hom}_{\mathcal{T}}(C,\mathcal{T}_{\leq -N})=0\,.$$
This implies that $\mathrm{Hom}_{\mathcal{T}}(C,X)=0$, for all compact
objects $C$ of $\mathcal{T}$. Since $\mathcal{T}$ is compactly
generated, we conclude that $X=0$. This proves the converse.
\end{proof}

\begin{lemma}\label{coproducts}
Let $(Y_p)_{p \in \mathbb{Z}}$ be in $\mathcal{T}$. We dispose of the
following isomorphism
$$ \mathrm{H}^n\,(\underset{p}{\coprod}Y_p)
\stackrel{\sim}{\longleftarrow}\underset{p}{\coprod}\, \mathrm{H}^n(Y_p)\,,
$$
for all $n \in \mathbb{Z}$.
\end{lemma}

\begin{proof}
By definiton $\mathrm{H}^n:= \tau_{\geq n}\,\tau_{\leq n}\,, n \in
\mathbb{Z}$. Since $\tau_{\geq n}$ admits a right adjoint, it is enough to show that $\tau_{\leq n}$
commute with infinite sums. 
We consider the following triangle
$$ \underset{p}{\coprod} \,\tau_{\leq n} Y_p \rightarrow
\underset{p}{\coprod}\,Y_p \rightarrow \underset{p}{\coprod}\, \tau_{>n} Y_p
\rightarrow S(\underset{p}{\coprod}\, \tau_{\leq n} Y_p)\,.$$
Here $\underset{p}{\coprod} \,\tau_{\leq n} Y_p$ belongs to $\mathcal{T}_{\leq
  n}$ since $\mathcal{T}_{\leq n}$ is stable under infinite sums. Let
$P$ be an object of $S^n\mathcal{U}$. Since $P$ is compact, we have
$$ \mathrm{Hom}_{\mathcal{T}}(P,\underset{p}{\coprod} \,\tau_{>n} Y_p)
\stackrel{\sim}{\longleftarrow} \underset{p}{\coprod}\,
\mathrm{Hom}_{\mathcal{T}}(P, \tau_{>n}Y_p)=0\,.$$
Since $\mathcal{T}_{\leq n}$ is generated by $S^n\mathcal{U}$,
$\underset{i}{\coprod} \,\tau_{>n} Y_p$ belongs to
$\mathcal{T}_{>n}$. Since the truncation triangle of
$\underset{p}{\coprod}\, Y_p$ is unique, this implies the following
isomorphism
$$ \underset{p}{\coprod} \,\tau_{\leq n}Y_p \stackrel{\sim}{\longrightarrow}
\tau_{\leq n}(\underset{p}{\coprod}\,Y_p)\,.$$
This proves the lemma.
\end{proof}

\begin{proposition}\label{filtration}
Let $X$ be an object of $\mathcal{T}$. Suppose that we are in the
conditions of proposition~\ref{extension} b). We dispose of the following isomorphism
$$ \underset{i}{\mathrm{hocolim}}\,\tau_{\leq i}X
\stackrel{\sim}{\longrightarrow} X\,.$$
\end{proposition}

\begin{proof}
We need only show that
$$\mathrm{H}^n(\underset{i}{\mbox{hocolim}}\, \tau_{\leq i}X)
\stackrel{\sim}{\longrightarrow} \mathrm{H}^n(X)\,,$$
for all $n \in \mathbb{Z}$. 
We dispose of the following triangle, \emph{cf.}~\cite{triangular},
$$ \underset{p}{\coprod}\,\tau_{\leq p}X \rightarrow \underset{q}{\coprod}\,
\tau_{\leq q}X \rightarrow \underset{i}{\mbox{hocolim}}\, \tau_{\leq
  i} X \rightarrow S(\underset{p}{\coprod}\,\tau_{\leq p}X)\,.$$
Since the functor $\mathrm{H}^n$ is homological, for all $n \in
\mathbb{Z}$ and it commutes with infinite sums by
lemma~\ref{coproducts}, we obtain a long exact sequence
$$ \cdots \rightarrow \underset{p}{\coprod}\, \mathrm{H}^n(\tau_{\leq
  p}X) \rightarrow \underset{q}{\coprod}\, \mathrm{H}^n(\tau_{\leq q}X)
  \rightarrow \mathrm{H}^n\, (\underset{i}{\mbox{hocolim}}\,
  \tau_{\leq i}X) \rightarrow \underset{p}{\coprod}\, \mathrm{H}^n\,
  S(\tau_{\leq p}X) \rightarrow \underset{q}{\coprod}\, \mathrm{H}^n\,
  S(\tau_{\leq q}X) \rightarrow \cdots$$
We remark that the morphism $\underset{p}{\coprod}\, \mathrm{H}^n\,
  S(\tau_{\leq p}X) \rightarrow \underset{q}{\coprod}\, \mathrm{H}^n\,
  S(\tau_{\leq q}X)$ is a split monomorphism and so we obtain $$ \mathrm{H}^n(X) = \underset{i}{\mbox{colim}}\,\mathrm{H}^n(\tau_{\leq i} X)
\stackrel{\sim}{\longrightarrow}
\mathrm{H}^n(\underset{i}{\mbox{hocolim}}\, \tau_{\leq i} X)\,.$$ 
\end{proof}

\end{document}